\documentclass[11pt]{amsart}

\usepackage{amsmath,amssymb,amsthm}
\usepackage{amsmath,amssymb,amsthm,amscd}
\usepackage[frame,cmtip,arrow,matrix,line,graph,curve]{xy}
\usepackage{graphpap, color}

\def\dfrac{\displaystyle\frac}
\def\dsum{\displaystyle\sum}

\newtheorem{prop}{Proposition}
\newtheorem{theo}[prop]{Theorem}
\newtheorem{lemm}[prop]{Lemma}
\newtheorem{coro}[prop]{Corollary}
\newtheorem{rema}[prop]{Remark}

\newcommand{\pa}{\partial}
\newcommand{\al}{\alpha}

\renewcommand{\leq}{\leqslant}
\renewcommand{\geq}{\geqslant}

\newcommand{\la}{\lambda}
\newcommand{\bs}{\boldsymbol}
\numberwithin{equation}{section}

\title{Newton-Maclaurin type inequalities for linear combinations of elementary symmetric functions}

\begin{document}
\author{Shuqi Hu}
\address{School of Mathematical Science\\
Jilin University\\ Changchun\\ China}
\email{husq21@mails.jlu.edu.cn}
\author{Changyu Ren}
\address{School of Mathematical Science\\ Jilin University\\ Changchun\\ China}
\email{rency@jlu.edu.cn}
\author{Ziyi Wang}
\address{School of Mathematical Science\\
Jilin University\\ Changchun\\ China}
\email{ziyiw23@mails.jlu.edu.cn}
\thanks{Research of the second author is supported by NSFC Grant No. 11871243.}

\begin{abstract}
In this paper, we establish Newton-Maclaurin type inequalities for
functions arising from linear combinations of primitively symmetric
polynomials. This generalization extends the classical
Newton-Maclaurin inequality to a broader class of functions.
\end{abstract}

 \maketitle
\vskip -6mm \par \qquad\hskip 1mm {\footnotesize Key Word: primitive
symmetric function; Newton inequality; Maclaurin inequa-
\par \qquad\hskip 2mm
lity; Hessian equation.
\par \qquad\hskip 2mm
MR2000: 05E05, 35A23.  }

\section{introduction}

In this paper, our primary focus is on Newton-Maclaurin type
inequalities pertaining to linear combinations of elementary
symmetric functions. First, define
$$\sigma_k({\bs x})={\dsum_{1\leq
i_{1}<\cdots<i_{k}\leq {n}}}x_ {i_{1}}\cdots x_ {i_{k}}$$ which is a
~$k$-th primitive symmetric function and the variables ${\bs
x}=(x_1,x_2,\cdots,x_{n})\in \mathbb{R}^n$. For example, when $n=3$,
$$\sigma_{1}({\bs x})=x_1+x_2+x_3,\quad \sigma_2({\bs
x})=x_1x_2+x_1x_3+x_2x_3,\quad \sigma_3({\bs x})=x_1x_2x_3.$$ It is
convenient to define $\sigma_0({\bs x})=1$, and define
$\sigma_k({\bs x})=0$ if $k<0$ or $k>n$. Furthermore, define a
 $k$-th primitive symmetric mean as
 $$E_k({\bs x})=\dfrac{\sigma_k({\bs x})}{C_{n}^k},\quad k=1,2,\cdots,{n},$$
 where ~$C_{n}^k=\dfrac{n!}{k!(n-k)!}$.

\par
 The classical Newton\cite{Newt} - Maclaurin\cite{Mac}
inequalities are
\begin{equation}\label{e1.1}
E_k({\bs x})^2\geq E_{k-1}({\bs x})E_{k+1}({\bs x}),\quad k=1,2\cdots,n,\quad \forall \bs x\in\mathbb{R}^n,
\end{equation}
\begin{equation}\label{e1.2}
E_1({\bs x})\geq E_2({\bs x})^{1/2}\geq\cdots\geq E_k({\bs x})^{1/k},\quad \forall\bs x\in\Gamma_{k},
\end{equation}
where $\Gamma_{k}$ is the G\r{a}rding's cone
$$\Gamma_{k}=\{x\in\mathbb{R}^n|\quad \sigma_{m}>0,\quad
m=1,2,\cdots,k\}.$$ All of the above inequalities are strict unless
$x_1=\cdots=x_n$, or both sides are $0$ at the same time. Over the
years, numerous scholars have contributed numerous generalizations
of Newton-Maclaurin type inequalities, including those presented in
\cite{RH,Men,NCP,SR}.

The Newton-Maclaurin inequalities play important roles in deriving
theoretical result for fully nonlinear partial differential
equations and geometric analysis. There are many important results
that need to use the Newton-Maclaurin inequalities, such as
\cite{BP,PX,PXF,GC}, etc. This is because the following $k$-Hessian
equations and curvature equations
\begin{align}
\sigma_k(\lambda(u_{ij}))&=f(x,u,\nabla u), \label{e1.03}\\
\sigma_k(\kappa(X))&=f(X,\nu)\label{e1.04}
\end{align}
are central studies in the field of fully nonlinear partial
differential equations and geometric analysis. Their left-hand side
$k$-Hessian operator $\sigma_k$ is the primary symmetric function
about the eigenvalues $\lambda=(\lambda_1,\cdots,\lambda_n)$ of the
Hessian matrix $(u_{ij})$ or
 the principal curvature
$\kappa=(\kappa_1,\cdots,\kappa_n)$ of the surface.

In recent years, the fully non-linear equations derived from linear
combinations of primary symmetric functions have received increasing
attention. For example, the following special Lagrangian equation
$${\rm Imdet}(\delta_{ij}+{\rm i}u_{ij})=\dsum_{k=0}^{[(n-1)/2]}(-1)^{k}\sigma_{2k+1}
(\lambda(u_{ij}))=0$$ was derived by Harvey and Lawson in their
study of the minimal submanifold problem\cite{HL}. The fully
nonlinear partial differential equations
$$P_m(u_{ij})=\sum_{k=0}^{m-1}(l_k^+)^{m-k}(x)\sigma_k(u_{ij})=g^{m-1}(x)
$$
were studied by Krylov \cite{K} and Dong \cite{Dong}. The following
curvature equations and complex Hessian equations
$$\dsum_{s=0}^k\alpha_s\sigma_k(\kappa(\chi))=f(\chi,\upsilon(\chi)),\quad
\chi\in M,$$
$$\dsum_{s=0}^k\alpha_s\sigma_k(\lambda)=f(z,u,Du)$$
were studied by Li-Ren-Wang \cite{LRW2} and Dong \cite{Dws}. For the
following curvature equations
$$\sigma_k(W_u(x))+\alpha
\sigma_{k-1}(W_u(x))=\dsum_{l=0}^{k-2}\alpha_l(x)\sigma_l(W_u(x)),\quad
x\in\mathbb{S}^n,$$ Guan-Zhang \cite{GZ} and Zhou \cite{Zhou}
studied the curvature estimates and the interior gradient estimates
respectively. Recently, Liu and Ren \cite{LR} discussed the
Pogorelov-type $C^{2}$ estimation for $(k-1)$-convex and $k$-convex
solutions of the following Sum Hessian equations
$$ \sigma_k(\lambda(u_{ij})) + \alpha\sigma_{k-1}(\lambda(u_{ij})) =
f(x, u,\nabla u).$$

The Newton-Maclaurin inequalities can be perceived as an intrinsic
characteristic of the $k$-Hessian operator. If the operators
situated on the left-hand side of the aforementioned equations can
establish a corresponding Newton-Maclaurin style inequality, it will
naturally facilitate further research on the equation. A natural
question is whether the Newton-Maclaurin type inequalities still
hold for the operators of linear combinations of these primary
symmetric functions? In a recent paper, we obtained the following
Newton-Maclaurin type inequalities for operators $E_k+\al E_{k-1}$
in \cite{Ren}:

\par
For any real number $\alpha\in\mathbb{R}$ and any
$\bs{x}=(x_{1},x_{2},\cdots,x_{n})\in\mathbb{R}^{n},$
 \begin{equation*}
 [E_{k}(\bs x)+\alpha E_{k-1}(\bs x)]^{2}\geq[E_{k-1}(\bs x)+\alpha E_{k-2}(\bs x)][E_{k+1}(\bs x)+\alpha E_{k}(\bs x)],\quad k=2\cdots,n-1.
 \end{equation*}
 Further, if $\alpha\geq0$ and $$E_{j}(\bs x)+\alpha E_{j-1}(\bs x)\geq0, \quad j=1\cdots,k,$$ then
$$[E_{1}(\bs x)+\alpha ]\geq [E_{2}(\bs x)+\alpha E_{1}(\bs x)]^{\frac{1}{2}}\geq\cdots\geq[E_{k}(\bs x)+\alpha E_{k-1}(\bs x)]^{\frac{1}{k}}.$$

A natural question is that if the Newton-Maclaurin type inequality
still holds for the functions which obtained by linear combination
of general primary symmestric functions. That is, for
$\al=(\la_1,\cdots,\al_n)\in \mathbb{R}^n, x=(x_1,\cdots,x_{n+1})\in
\mathbb{R}^{n+1}$, whether the following inequalities hold?
\begin{equation}\label{e1.05}
\left(\dsum_{k=1}^n \al_kE_k(x)\right)^2\geq\left(\dsum_{k=1}^n
\al_k E_{k-1}(x)\right)\left(\dsum_{k=1}^n\al_k E_{k+1}(x)\right).
\end{equation}

A counterexample was given in \cite{Ren} to illustrate that
inequality \eqref{e1.05} generally does not hold. So can we add some
conditions to $\al$ to ensure inequality \eqref{e1.05} hold?

In this paper, we consider the Newton-Maclaurin type inequalities
for operators of three elementary symmetric functions $[E_{k}(\bs
{x})+a E_{k-1}(\bs x)+b E_{k-2}(\bs x)]$. For convenience, we define
$$S_{k}(\bs x):=E_{k}(\bs x)+(\alpha+\beta)E_{k-1}(\bs
x)+\alpha\beta E_{k-2}(\bs x),$$ where $\alpha, \beta\in\mathbb{R}$
and $\bs x=(x_{1},x_{2},\cdots,x_{n})\in\mathbb{R}^{n}$.

\begin{theo}\label{th1}
For any real $\alpha, \beta\in\mathbb{R}$ and any $\bs
x\in\mathbb{R}^{n}$, we have
 \begin{align}\label{e1.3}
 S_{k}^{2}(\bs x)\geq S_{k-1}(\bs x)S_{k+1}(\bs x),\quad
 k=3,4,\cdots,n-1.
 \end{align}
 The inequality is strict unless $x_{1}=x_{2}=\cdots=x_{n},$
or unless both sides of the inequality are $0$,  or unless
 \begin{equation}\label{e1.4}
 \frac{S_{k}(\bs x)}{S_{k-1}(\bs x)}=\frac{S_{k+1}(\bs x)}{S_{k}(\bs x)}=-\al=-\beta.
 \end{equation}

 \par
 Further, if $\alpha,\beta \geq 0, E_1(\bs x)\geq 0, E_2(\bs x)\geq 0$, and
  \begin{equation*}
 S_{m}(\bs x)\geq 0,\quad {\rm for~~ all}~~ m=3,4\cdots,k.
\end{equation*}
 Then
\begin{equation}\label{e1.5}
S_1(\bs x) \geq S_{2}^{1/2}(\bs x) \geq \cdots\geq S_{k}^{1/k}(\bs
x),\quad k=2,3,\cdots,n.
\end{equation}
\end{theo}

\par
\begin{rema}\label{r2}
There is a special case when the equation \eqref{e1.4} holds: when
$\alpha=\beta$ and there are $n-2$ elements in $x_1,x_2,\cdots,x_n$
take the value $-\alpha$.
\end{rema}

Obviously, Theorem $\ref{th1}$ is a generalization of classical
Newton-Maclaurin inequality and Newton-Maclaurin  type inequality in
\cite{Ren}. The structural condition to constants $a$ and $b$ in
operators $[E_{k}(\bs {x})+a E_{k-1}(\bs x)+b E_{k-2}(\bs x)]$ can
be seen as the following polynomial of $t$
\begin{equation}\label{e1.6}
 t^2+a
t+b=0
\end{equation}
only has real roots. And its real roots are $-\al$ and $-\beta$. The
polynomial \eqref{e1.6} which only has real roots is a sufficient
condition of Theorem \ref{th1}. Probably, this condition is also
necessary. In Second 2, we will provide a counterexample for the
case of $k=3$. For the case of $k>3$, we can only provide partial
counterexamples.

\par
Similar as \cite{G} for $E_{k}(\bs x)$, from \eqref{e1.3} we have
\begin{coro}\label{cor4}
Let $3\leq l<  k\leq n$, assume
 \begin{equation}\nonumber
 S_{q}(\bs x)\geq 0,\quad\text{for all}~q=l,\cdots,k-1,
  \end{equation}
then
\begin{align*}
&S_{l}(\bs x)S_{k-1}(\bs x) \geq S_{l-1}(\bs x)S_{k}(\bs x).
\end{align*}
\end{coro}

\par

For the more terms case, we can establish the Newton-Maclaurin type
inequalities for operators
\begin{equation}\label{e1.7}
S_{k;s}(\bs x):=C^{0}_{s}E_{k}(\bs x)+C^{1}_{s}\alpha E_{k-1}(\bs
x)+\cdots+C^{s}_{s}\alpha^s E_{k-s}(\bs x),
\end{equation}
where
$\alpha\in\mathbb{R}$ and $\bs
x=(x_{1},x_{2},\cdots,x_{n})\in\mathbb{R}^{n}$.

\begin{theo}\label{th5}
For general $s\geq1$, any real $\alpha\in\mathbb{R}$ and any $\bs
x\in\mathbb{R}^{n}$, we have
 \begin{align}\label{e1.8}
 S_{k;s}^{2}(\bs x)\geq S_{k-1;s}(\bs x)S_{k+1;s}(\bs x), \quad  k=s+1,\cdots, n-1.
 \end{align}
 The inequality is strict unless $x_{1}=x_{2}=\cdots=x_{n},$
or unless both sides of the inequality are $0$,  or unless
 \begin{equation}\label{e1.9}
 \frac{S_{k;s}(\bs x)}{S_{k-1;s}(\bs x)}=\frac{S_{k+1;s}(\bs x)}{S_{k;s}(\bs x)}=-\al.
 \end{equation}

 \par
 Further, if $\alpha \geq 0, E_1(\bs x)\geq 0, E_2(\bs x)\geq 0,\cdots,E_s(\bs x)\geq 0$, and
  \begin{equation*}
 S_{m;s}(\bs x)\geq 0,\quad {\rm for~~ all}~~ m=s+1,s+2,\cdots,k.
\end{equation*}
 Then \eqref{e1.8} hold for all $k=1,\cdots,n-1$, and
\begin{equation}\label{e1.10}
S_{1;s}(\bs x) \geq S_{2;s}^{1/2}(\bs x)\geq\cdots \geq
S_{k;s}^{1/k}(\bs x),\quad k=2,3,\cdots,n.
\end{equation}
\end{theo}

\begin{rema}
There is a special case when the equation \eqref{e1.9} holds: when
there are $n-s$ elements in $x_1,x_2,\cdots,x_n$ take the value
$-\alpha$.
\end{rema}

\par
Similar as \cite{G} for $E_{k}(\bs x)$, from \eqref{e1.8} we have
\begin{coro}\label{cor7}
Let $s+1 \leq l<  k\leq n$, assume
 \begin{equation}\nonumber
 S_{q;s}(\bs x)\geq 0,\quad\text{for all}~q=l,\cdots,k-1,
  \end{equation}
then
\begin{align*}
&S_{l;s}(\bs x)S_{k-1;s}(\bs x) \geq S_{l-1;s}(\bs x)S_{k;s}(\bs x).
\end{align*}
\end{coro}

The organization of this paper is as follows. In Section 2, we will
prove Theorem \ref{th1} and in Section 3 we will prove Theorem
\ref{th5}.

\section{The proof of Theorem \ref{th1}}
First, we consider the following inequality in Euclidean space
$\mathbb{R}^{4}$.

\par

\begin{lemm}\label{lem2.1}
For any real numbers $\alpha, \beta\in\mathbb{R}$ and $\bs
z=(z_{1},z_{2},z_{3},z_{4})\in\mathbb{R}^{4}$, we have
 \begin{align}\label{e2.1}
 S_{3}^{2}(\bs z)\geq S_{2}(\bs z)S_{4}(\bs z),
 \end{align}
 the inequality is strict unless both sides of the inequality sign are $0$, or unless $z_{1}=z_{2}=z_{3}=z_{4}$, or unless
  \begin{equation}\nonumber
  \frac{S_{3}(\bs z)}{S_{2}(\bs z)}=\frac{S_{4}(\bs z )}{S_{3}(\bs z)}=-\al=-\beta.
  \end{equation}
\end{lemm}

Proof: The following is one of the most straightforward and simplest calculations
\begin{align}
&576[S_{3}^{2}(\bs z)- S_{2}(\bs z)S_{4}(\bs z)]  \nonumber\\
=&576\Big\{[E_{3}(\bs z)+(\alpha+\beta )E_{2}(\bs z)+\alpha\beta E_{1}(\bs z)]^{2}  \nonumber\\
&-[E_{2}(\bs z)+(\alpha+\beta )E_{1}(\bs z)+\alpha\beta][E_{4}(\bs z)+(\alpha+\beta )E_{3}(\bs z)+\alpha\beta E_{2}(\bs z)]\Big\} \nonumber\\
=&3\Big\{[(z_{1}-z_{2})(z_{3}+\alpha)(z_{4}+\beta)+(z_{1}-z_{2})(z_{4}+\alpha)(z_{3}+\beta)]^{2}  \label{equa}\\
&+[(z_{1}-z_{3})(z_{2}+\alpha)(z_{4}+\beta)+(z_{1}-z_{3})(z_{4}+\alpha)(z_{2}+\beta)]^{2}  \nonumber\\
&+[(z_{1}-z_{4})(z_{2}+\alpha)(z_{3}+\beta)+(z_{1}-z_{4})(z_{3}+\alpha)(z_{2}+\beta)]^{2}  \nonumber\\
&+[(z_{2}-z_{3})(z_{1}+\alpha)(z_{4}+\beta)+(z_{2}-z_{3})(z_{4}+\alpha)(z_{1}+\beta)]^{2}  \nonumber\\
&+[(z_{2}-z_{4})(z_{1}+\alpha)(z_{3}+\beta)+(z_{2}-z_{4})(z_{3}+\alpha)(z_{1}+\beta)]^{2}  \nonumber\\
&+[(z_{3}-z_{4})(z_{1}+\alpha)(z_{2}+\beta)+(z_{3}-z_{4})(z_{2}+\alpha)(z_{1}+\beta)]^{2}\Big\}  \nonumber\\
&+2\Big\{[(z_{1}-z_{2})(z_{3}-z_{4})(\alpha-\beta)]^2  \nonumber\\
&+[(z_{1}-z_{3})(z_{2}-z_{4})(\alpha-\beta)]^2  \nonumber\\
&+[(z_{1}-z_{4})(z_{2}-z_{3})(\alpha-\beta)]^2\Big\}\geq0. \nonumber
\end{align}
It is obvious to see that the equality sign holds when
$z_{1}=z_{2}=z_{3}=z_{4}$, or $\alpha=\beta$ and any two elements of
$z_{1},z_{2},z_{3},z_{4}$ are equal to $-\alpha$. We assume that
$\alpha=\beta$ and ~$z_{1}=z_{2}=-\alpha$, then we get
\begin{align*}
&[E_{2}(\bs z)+(\alpha+\beta )E_{1}(\bs z)+\alpha\beta E_{0}(\bs z)]=\frac{1}{6}(z_{3}z_{4}+\alpha z_{3}+\alpha z_{4}+\alpha^{2}),\\
&[E_{3}(\bs z)+(\alpha+\beta )E_{2}(\bs z)+\alpha\beta E_{1}(\bs z)]=-\frac{1}{6}\alpha (z_{3}z_{4}+\alpha z_{3}+\alpha z_{4}+\alpha^{2}),\\
&[E_{4}(\bs z)+(\alpha+\beta )E_{3}(\bs z)+\alpha\beta E_{2}(\bs
z)]=\frac{1}{6}\alpha^{2}(z_{3}z_{4}+\alpha z_{3}+\alpha
z_{4}+\alpha^{2}).
\end{align*}
Thus we complete the proof of Lemma \ref{lem2.1}.

\begin{rema} The calculation of the equation \eqref{equa} is trivial but
tedious. The same goes for some equations on page 9. Therefore, we
suggest that readers use mathematical software (For example,
Mathematica) to assist in the checking.
\end{rema}

The following lemma is an useful tool to prove Newton's inequalities
from \cite{S,G,NCP}.

\begin{lemm}\label{lem2.2}
If
\begin{equation}\nonumber
F(x,y)=c_{0}x^{n}+c_{1}x^{n-1}y+\cdots+c_{n}y^{n},
\end{equation}
is a homogeneous function of the n-th degree in x and y whose roots
are all real for x and y, then the same is true for all
non-identical $0$ equations
\begin{equation}\nonumber
\frac{\pa^{i+j}F}{\pa x^{i}\pa y^{j}}=0,
\end{equation}
obtained from it by partial  differentiation with respect to $x$ and
$y$. Further, if  $Q$ is one of these equations, and it has a
multiple root $\gamma$, then $\gamma$ is also a higher multiplicity
root of the equation which is derived from the differentiation of
$Q$.
\end{lemm}
\par
Now we prove the Theorem \ref{th1}.

\par\noindent
 {\bf Proof of Theorem \ref{th1}}.
Assuming that $P(t)$ is an $n$-degree polynomial with real roots
$x_{1},x_{2},\cdots,x_{n}$. Then $P(t)$ can be represented as
 \begin{equation}\nonumber
 P(t)=\prod_{i=1}^{n}(t-x_{i})=E_{0}(\bs x)t^{n}-C_{n}^{1} E_{1}(\bs x) t^{n-1}+C_{n}^{2} E_{2}(\bs x) t^{n-2}-\cdots+(-1)^{n} E_{n}(\bs x),
 \end{equation}
 We apply the lemma \ref{lem2.2} to the following related homogeneous polynomials
\begin{equation}\nonumber
F(t,s)=E_{0}(\bs x)t^{n}-C_{n}^{1} E_{1}(\bs x) t^{n-1} s+C_{n}^{2}
E_{2}(\bs x) t^{n-2} s^{2}-\cdots+(-1)^{n} E_{n}(\bs x) s^{n} .
\end{equation}
\par

Consider the case of the derivative $\dfrac{\partial^{n-4}
F}{\partial t^{n-1-k} \partial s^{k-3}}$, where $k=3, \cdots, n-1$.
We obtain the following family of quartic polynomials
\begin{equation}\label{e2.2}
E_{k-3}(\bs x)t^{4}-4E_{k-2}(\bs x)t^{3}s+6E_{k-1}(\bs
x)t^{2}s^{2}-4E_{k}(\bs x)ts^{3}+E_{k+1}(\bs x)s^{4},
\end{equation}
where $k=3,\cdots,n-1$. The above polynomials only have real roots
by lemma \ref{lem2.2}.

Next we will divide it into two cases to deal with \eqref{e1.3}:
\begin{center}
$E_{k+1}(\bs x)\not=0$; \quad $E_{k+1}(\bs x)=0$.
\end{center}

{\bf Case~A: } When $E_{k+1}(\bs x)\not=0$, from \eqref{e2.2}, we
get the polynomial

\begin{equation}\nonumber
s^{4}-\frac{4 E_{k}({\bs x})}{E_{k+1}({\bs x})} s^{3}+\frac{6
E_{k-1}({\bs x})}{E_{k+1}({\bs x})}s^{2}-\frac{4E_{k-2}({\bs
x})}{E_{k+1}({\bs x})}s+\frac{E_{k-3}({\bs x})}{E_{k+1}({\bs x})},
\end{equation}
which has four real roots. We denote the roots by
${z}_{1},{z}_{2},{z}_{3},{z}_{4}$, and denote ${\bs
z}=({z}_{1},{z}_{2},{z}_{3},{z}_{4})$, then we have
\begin{equation}\nonumber
E_{1}({\bs z})=\frac{E_{k}({\bs x})}{E_{k+1}({\bs x})}, ~~E_{2}({\bs
z})=\frac{E_{k-1}({\bs x})}{E_{k+1}({\bs x})}, ~~ E_{3}({\bs z}
)=\frac{E_{k-2}({\bs x})}{E_{k+1}({\bs x})}, ~~ E_{4}({\bs
z})=\frac{E_{k-3}({\bs x})}{E_{k+1}({\bs x})}.
\end{equation}
Replacing $\alpha$ and $\beta$ with $\dfrac{1}{\alpha}$ and
$\dfrac{1}{\beta}$ respectively, and using the Lemma \ref{lem2.1} we
obtain \eqref{e1.3}.

{\bf Case~B:}
  When $E_{k+1}(\bs x)=0$, we first consider the case of the derivative $\dfrac{\partial^{n-1-k}
F}{\partial t^{n-1-k}}$. We obtain the following $k+1$-order
polynomial
\begin{equation}\label{e2.3}
E_{0}(\bs x)t^{k+1}-C_{k+1}^{1} E_{1}(\bs x) t^{k}+C_{k+2}^{2}
E_{2}(\bs x) t^{k-1}-\cdots+(-1)^{k+1} E_{k+1}(\bs x),
\end{equation}
and the polynomial only has real roots by lemma \ref{lem2.2}. We
denote the roots by ${y}_{1},{y}_{2},\cdots,{y}_{k+1}$, and denote
${\bs y}=({y}_{1},{y}_{2},\cdots,{y}_{k+1})$. Since
$$
0=E_{k+1}(\bs x)=E_{k+1}(\bs y)={y}_{1}\cdot{y}_{2}\cdots y_{k+1},
$$
this shows that the polynomial $\eqref{e2.3}$ has roots of value 0,
and we denote these 0 value roots as
${y}_{1},{y}_{2},\cdots,{y}_{r}$ respectively. Let
$\mathbf{e}=(1,\cdots,1,0,\cdots,0)$ represent a vector where the
first $r$ elements are 1 and the last $k+1-r$ elements are 0.
 Let $\widetilde{\bs y}=\bs y+\varepsilon\mathbf{e}=(\varepsilon,\cdots,\varepsilon,y_{r+1},\cdots,y_{k+1})$, then for any $ \varepsilon>0$, we have $E_{k+1}(\widetilde{\bs y})\not= 0,$ and
 $$
\lim_{\varepsilon\rightarrow0}E_{m}(\widetilde{\bs y})=E_{m}( \bs
y)=E_{m}( \bs x),\quad {\rm for~~all}~~ m=1,\cdots,k+1.
$$
In the following we consider the $k + 1$-order polynomial
 \begin{align*}
  E_{0}(\widetilde{\bs y})t^{k+1}-C_{k+1}^{1}E_{1}({\widetilde{\bs y}})t^{k}s+C_{k+2}^{2}E_{2}({\widetilde{\bs y}})t^{k-1}s^{2}-\cdots+(-1)^{k+1}E_{k+1}({\widetilde{\bs y}})s^{k+1}.
  \end{align*}
 By the definition of $\widetilde{\bs y}$, we know that this polynomial only has real roots. Differentiating the polynomial above with respect to $s$ $k-3$ times, we obtain
   \begin{align}\label{e2.4}
  E_{k-3}(\widetilde{\bs y})t^{4}-4E_{k-2}({\widetilde{\bs y}})t^{3}s+6E_{k-1}({\widetilde{\bs y}})t^{2}s^{2}-4E_{k}({\widetilde{\bs y}})ts^{3}+E_{k+1}({\widetilde{\bs y}})s^{4},
  \end{align}
 where $k=3,\cdots,n-1$.

Since $E_{k+1}(\widetilde{\bs y})\not= 0,$ by \eqref{e2.4}, we have
the polynomial \begin{equation}\nonumber s^{4}-\frac{4
E_{k}(\widetilde{\bs y})}{E_{k+1}(\widetilde{\bs y})} s^{3}+\frac{6
E_{k-1}(\widetilde{\bs y})}{E_{k+1}(\widetilde{\bs y})}
s^{2}-\frac{4E_{k-2}(\widetilde{\bs y})}{E_{k+1}(\widetilde{\bs
y})}s+\frac{E_{k-3}(\widetilde{\bs y})}{E_{k+1}(\widetilde{\bs y})},
\end{equation} which only has real roots.

We denote the roots by
$\widetilde{z}_{1},\widetilde{z}_{2},\widetilde{z}_{3},\widetilde{z}_{4}$
and denote $\widetilde{\bs z}=(\widetilde{z}_{1},\widetilde{z}
_{2},\widetilde{z}_{3},\widetilde{z}_{4})$, then
\begin{equation}\nonumber
E_{1}(\widetilde{\bs z})=\frac{E_{k}(\widetilde{\bs
y})}{E_{k+1}(\widetilde{\bs y})}, \quad E_{2}(\widetilde{\bs
z})=\frac{E_{k-1}(\widetilde{\bs y})}{E_{k+1}(\widetilde{\bs y})},
\quad E_{3}(\widetilde{\bs z})=\frac{E_{k-2}(\widetilde{\bs
y})}{E_{k+1}(\widetilde{\bs y})}, \quad E_{4}(\widetilde{\bs
z})=\frac{E_{k-3}(\widetilde{\bs y})}{E_{k+1}(\widetilde{\bs y})}.
\end{equation}
Using the lemma \ref{lem2.1} we get
 \begin{align*}
 S_{k}^{2}(\widetilde{\bs y})
 \geq S_{k-1}(\widetilde{\bs y})S_{k+1}(\widetilde{\bs y}),
 \end{align*}
Taking $\varepsilon\rightarrow0$, we have $\widetilde{\bs
y}\rightarrow \bs y$,  which gives
$$E_{i}(\widetilde{\bs y})\rightarrow
E_{i}(\bs y)=E_{i}(\bs x), ~~ i=1,2,\cdots,k+1,$$ then we obtain
 \begin{align*}
 S_{k}^{2}(\bs x)
 \geq S_{k-1}(\bs x)S_{k+1}(\bs x).
 \end{align*}

\par
By Lemma $\ref{lem2.1}$, the inequalities of \eqref{e1.3} are strict
unless $x_{1}=x_{2}=\cdots=x_{n}$, or unless
 $$\frac{S_{k}(\bs x)}{S_{k-1}(\bs x)}=\frac{S_{k+1}(\bs x)}{S_{k}(\bs x)}=-\al=-\beta.$$
A special case in which the above equality holds is that
$\alpha=\beta$ and there are $n-2$ elements of $x_1,x_2,\cdots,x_n$
valued $-\alpha$.

\par
Now, we prove \eqref{e1.5}. Notice that $E_0(\bs x)=1, E_k(\bs x)=0$
for $k<0$, when $\al,\beta\geq 0,E_1(\bs x)\geq 0,E_2(\bs x)\geq 0$,
we have
\begin{align}\label{e2.5}
&S_{1}^2(\bs x)-S_{0}(\bs x)S_{2}(\bs x)\\
=&[E_1+(\alpha+\beta)]^2-[E_2+(\alpha+\beta)E_1+\alpha\beta] \nonumber\\
 =&E_1^2-E_2+(\alpha+\beta)E_1+(\alpha+\beta)^2-\alpha\beta\geq 0, \nonumber
\end{align}
and
\begin{align*}
&S^2_{2}(\bs x)-S_{1}(\bs x)S_{3}(\bs x)\\
=&E_2^2+(\al^2+\beta^2+\al\beta)E_1^2+\al^2\beta^2+(\al+\beta)E_1E_2+\al\beta(\al+\beta)E_1\\
&-[E_1E_3+(\al+\beta)E_3+(\al^2+\beta^2)E_2]\\
\geq&\al\beta E_1^2+\al^2\beta^2+\alpha\beta(\alpha+\beta)E_1\geq 0.
\end{align*}
\par
 So the proof of \eqref{e1.5} is similar as \cite{G,Ren}.  \eqref{e2.5}
 imply that
 $$S_{1}(\bs x)\geq S_{2}^{1/2}(\bs x).$$
For $2\leq k\leq n-1$, we assume
$$S_{k-1}^{1/{(k-1)}}(\bs x) \geq S_{k}^{1/k}(\bs x).$$
Combining the inequality \eqref{e1.3} and the inequalities above, we
get
\begin{align*}
S_{k}^2(\bs x) \geq S_{k-1}(\bs x)S_{k+1}(\bs x) \geq
S_{k}^{(k-1)/k}(\bs x)S_{k+1}(\bs x).
\end{align*}
Thus,
\begin{equation*}
S_{k}^{1/k}(\bs x)\geq S_{k+1}^{1/(k+1)}(\bs x).
\end{equation*}
The inequality $\eqref{e1.5}$ is proved. \hspace*{\fill}$\Box$

\vskip 3mm\par

{\bf $\bullet$ Counterexample}
\par
Assuming polynomial $t^2+at+b$ has two complex roots $-c\pm  {\rm
i}d$ with $d\neq 0$, then
$$S_k'(x)=[E_{k}(\bs {x})+a E_{k-1}(\bs
x)+b E_{k-2}(\bs x)]$$
can be rewrite as
\begin{equation*}
S_k'(x)=E_{k}(\bs {x})+2c E_{k-1}(\bs x)+(c^2+d^2) E_{k-2}(\bs x).
\end{equation*}
Next, we will provide a counterexample for the case of $k=3$. That
is to say, we found a $\bs z_0\in \mathbb{R}^n$ which satisfies
$$
[S_3'(\bs z_0)]^2<S_{2}'(\bs z_0)S_{4}'(\bs z_0).
$$
\par
Denote $\bs z'=(z_3,\cdots,z_n)$, by (iii) of section 2 in
\cite{RW2}, we have
\begin{equation}\label{e2.6}
\sigma_k(\bs z)=z_1z_2\sigma_{k-2}(\bs z')+(z_1+z_2)\sigma_{k-1}(\bs
z')+\sigma_{k}(\bs z').
\end{equation}
We divide it into the following cases.
\par
{\bf Case 1.} $|c|\geq |d|$. We choose $\bs
z_0=((n-3)c,(n-3)c,-c,-c,\cdots,-c)$, using \eqref{e2.6}, we have
$$E_1(\bs z_0)=\dfrac{2(n-3)c-(n-2)c}{n}=\dfrac{(n-4)c}{n},$$
\begin{align*}
E_2(\bs
z_0)=&\dfrac{c^2}{C_n^2}[(n-3)^2-2(n-3)(n-2)+C_{n-2}^2]=\dfrac{-(n-3)c^2}{n-1},\\
E_3(\bs
z_0)=&\dfrac{c^3}{C_n^3}[-(n-k)^2(n-2)+2(n-k)C_{n-2}^2-C_{n-2}^3]=\dfrac{-(n-3)(n-4)c^3}{n(n-1)},\\
E_4(\bs
z_0)=&\dfrac{c^4}{C_n^4}[(n-k)^2C_{n-2}^2-2(n-k)C_{n-2}^3+C_{n-2}^4]=\dfrac{(5n^2-25n+32)c^4}{n(n-1)},
\end{align*}
and
\begin{align*}
S_2'(\bs z_0)=&E_2(\bs z_0)+2c E_1(\bs
z_0)+(c^2+d^2)=\dfrac{2(n-2)^2c^2+n(n-1)d^2}{n(n-1)},\\
S_3'(\bs z_0)=&\dfrac{-2(n-2)^2c^3+(n^2-5n+4)cd^2}{n(n-1)},\\
S_4'(\bs z_0)=&\dfrac{2(n-2)^2c^4-n(n-3)c^2d^2}{n(n-1)}.
\end{align*}
 So we have
$$
[S_3'(\bs z_0)]^2-S_{2}'(\bs z_0)S_{4}'(\bs
z_0)=\dfrac{2(n-2)^3c^2d^2[-2(n-2)c^2+(n-1)d^2]}{n^2(n-1)^2}<0.
$$

\par
{\bf Case 2.} $|c|<|d|$. We choose $\bs
z_0=((n-3)d,(n-3)d,-d,-d,\cdots,-d)$, using \eqref{e2.6}, we have
$$
\begin{array}{ll} E_1(\bs z_0)=\dfrac{(n-4)d}{n}, &
E_2(\bs
z_0)=\dfrac{-(n-3)d^2}{n-1},\\
E_3(\bs z_0)=\dfrac{-(n-3)(n-4)d^3}{n(n-1)},\quad  & E_4(\bs
z_0)=\dfrac{(5n^2-25n+32)d^4}{n(n-1)},
\end{array}
$$
and
\begin{align*}
S_2'(\bs z_0)=&\dfrac{n(n-1)c^2+2(n^2-5n+4)cd+2nd^2}{n(n-1)},\\
S_3'(\bs z_0)=&\dfrac{(n^2-5n+4)c^2d-2n(n-3)cd^2+2(n-4)d^3}{n(n-1)},\\
S_4'(\bs
z_0)=&\dfrac{-n(n-3)c^2d^2-2(n^2-7n+12)cd^3+2(2n^2-11n+16)d^4}{n(n-1)}.
\end{align*}
So we have
\begin{align*}
&[S_3'(\bs z_0)]^2-S_{2}'(\bs z_0)S_{4}'(\bs
z_0)\\
=&\dfrac{2(n-2)^3d^2[(n-1)c^4+2(n-4)c^2d^2-4(n-4)cd^3-4d^4]}{n^2(n-1)^2}\\
\leq&\dfrac{2(n-2)^3d^2[(n-1)c^4-2(n-4)cd^3-4d^4]}{n^2(n-1)^2}<0.
\end{align*}

\par
\begin{rema}
For the case of $k=4$, when $|c|\geq |d|$, we can choose
\begin{equation}\label{e2.8}
\bs z_0=((n-k)c,(n-k)c,-c,-c,\cdots,-c),\end{equation} which
satisfies
\begin{equation}\label{e2.9}
[S_k'(\bs z_0)]^2-S_{k-1}'(\bs z_0)S_{k+1}'(\bs
z_0)<0.\end{equation} For the case of $k=5$, when $|c|\geq |d|$ and
$n\leq 9$, the $\bs z_0$ in \eqref{e2.8} also satisfies inequalities
\eqref{e2.9}.
 Since the calculations similar as $k=3$, here
we omitted the calculation process. We suggest that readers verify
it by using mathematical software.
\end{rema}

\vskip 3mm

\section{The proof of Theorem \ref{th5}}

\par
To prove the Theorem $\ref{th5}$, we need to prove the following
lemma.
\begin{lemm}\label{lem3.1}
For general $k\geq2,~s=k-1,~n=k+1$, any real $\alpha\in\mathbb{R}$
and any $\bs z\in\mathbb{R}^{n}$, we have
\begin{equation}\label{e3.1}
S_{k;s}^{2}(\bs z)\geq  S_{k-1;s}(\bs z)S_{k+1;s}(\bs z),
\end{equation}
the inequality is strict unless both sides of the inequality sign
are $0$, or unless $z_{1}=z_{2}=\cdots=z_{n}$, or unless
\begin{equation}\nonumber
  \frac{S_{k;s}(\bs z)}{S_{k-1;s}(\bs z)}=\frac{S_{k+1;s}(\bs z )}{S_{k;s}(\bs z)}=-\al.
\end{equation}
\end{lemm}

{\bf Proof}. Let's first prove two equations.
\par
{\bf Equation 1.} For $n=k+1$,
\begin{equation}\label{e3.2}
kn^2(E_k^2-E_{k-1}E_{k+1})(\bs z)=\sum_{1\leq i<j\leq
n}(z_i-z_j)^2\prod_{l\neq i,j}z_l^2.
\end{equation}
We can obtain this equation by directly calculation. For
convenience, we can always rewrite
$$\sigma_{n-1}=\dsum_{i=1}^n\frac{\sigma_n}{z_i},\quad \sigma_{n-2}=\dsum_{1\leq i<j\leq n}\frac{\sigma_n}{z_iz_j}$$
regardless of whether $z_i$ and $z_j$ are 0 or not.  Note that
$k=n-1$, the left hand side of \eqref{e3.2} can be represented as
\begin{align*}
kn^2(E_k^2-E_{k-1}E_{k+1})(\bs z)
=&(n-1)\sigma_{n-1}^2-2n\sigma_{n-2}\sigma_n\\
=&(n-1)\Big(\sum_{i=1}^n\frac{\sigma_n}{z_i}\Big)^2-2n\sum_{1\leq i<j\leq n}\frac{\sigma_n}{z_iz_j}\sigma_n\\
=&\sigma_n^2\Big[(n-1)\Big(\sum_{i=1}^n\frac{1}{z_i}\Big)^2-2n\sum_{1\leq i<j\leq n}\frac{1}{z_iz_j}\Big]\\
=&\sigma_n^2\Big[(n-1)\Big(\sum_{i=1}^n\frac{1}{z_i^2}\Big)-2\sum_{1\leq
i<j\leq n}\frac{1}{z_iz_j}\Big].
\end{align*}
Meanwhile, the right hand side of \eqref{e3.2}
\begin{align*}
\sum_{1\leq i<j\leq n}(z_i-z_j)^2\prod_{l\neq i,j}z_l^2
=&\sum_{1\leq i<j\leq n}(z_i-z_j)^2\frac{\sigma_n^2}{z_i^2z_j^2}\\
=&\sigma_n^2\sum_{1\leq i<j\leq n}(\frac{1}{z_i}-\frac{1}{z_j})^2\\
=&\sigma_n^2\Big[(n-1)\Big(\sum_{i=1}^n\frac{1}{z_i^2}\Big)-2\sum_{1\leq
i<j\leq n}\frac{1}{z_iz_j}\Big].
\end{align*}
Comparing the left hand side with the right hand side, we obtain
\eqref{e3.2}.

\par
{\bf Equation 2.}  For $n=k+1, s=k-1$,
\begin{equation}\label{e3.3}
(S_{k;s}^2-S_{k-1;s}S_{k+1;s})(\bs z)=(E_k^2-E_{k-1}E_{k+1})(\bs
z+\alpha\bs e).
\end{equation}
We give the expansions of each term in the equation above. By the
define of $S_{k;s}$ in \eqref{e1.7},
\begin{align*}
S_{k;s}(\bs z)&=\sum_{i=k-s}^kC^{i-1}_{s}\alpha^{k-i}E_i(\bs z)=\sum_{i=1}^k\frac{i(k+1-i)}{(k+1)k}\alpha^{k-i}\sigma_i,\\
S_{k-1;s}(\bs z)&=\sum_{i=k-1-s}^{k-1}C^{i}_{s}\alpha^{k-1-i}E_i(\bs z)=\sum_{i=0}^{k-1}\frac{(k-i)(k+1-i)}{(k+1)k}\alpha^{k-1-i}\sigma_i,\\
S_{k+1;s}(\bs
z)&=\sum_{i=k+1-s}^{k+1}C^{i-2}_{s}\alpha^{k+1-i}E_i(\bs
z)=\sum_{i=2}^{k+1}\frac{i(i-1)}{(k+1)k}\alpha^{k+1-i}\sigma_i.
\end{align*}
Using the well known result
$$
E_k(\bs z+\alpha\bs e)=\sum_{i=0}^{k} C^{i}_{k}\alpha^{k-i}E_i(\bs
z),
$$
where $\bs e=(1,1,\cdots,1)\in \mathbb{R}^n$, see \cite{Gard}, we
have
\begin{align*}
E_k(\bs z+\alpha\bs e)&=\sum_{i=0}^{k} C^{i}_{k}\alpha^{k-i}E_i(\bs z)=\sum_{i=0}^k\frac{k+1-i}{k+1}\alpha^{k-i}\sigma_i,\\
E_{k-1}(\bs z+\alpha\bs
e)&=\sum_{i=0}^{k-1}C^{i}_{k-1}\alpha^{k-1-i}E_i(\bs
z)=\sum_{i=0}^{k-1}\frac{(k-i)(k+1-i)}{(k+1)k}\alpha^{k-1-i}\sigma_i,\\
E_{k+1}(\bs z+\alpha\bs e)&=\sum_{i=0}^{k+1}
C^{i}_{k+1}\alpha^{k+1-i}E_i(\bs
z)=\sum_{i=0}^{k+1}\alpha^{k+1-i}\sigma_i.
\end{align*}

From the expansion above, we only need to prove that all the terms
containing $\alpha$ at both hand sides of \eqref{e3.3} are equal.
\par
On the left hand side of \eqref{e3.3}, the coefficient of the term
$\alpha^{2k-b}(b\geq2)$ is
\begin{align*}
L_{2k-b}&=\sum_{\substack{i+j=b\\1\leq i\leq k\\1\leq j\leq
k}}\frac{(k+1-i)(k+1-j)ij}{(k+1)^2k^2}\sigma_i\sigma_j-\sum_{\substack{p+q=b\\0\leq
p\leq k-1\\2\leq q\leq
k+1}}\frac{(k+1-p)(k-p)q(q-1)}{(k+1)^2k^2}\sigma_p\sigma_q\\
&:=L_{2k-b}^1-L_{2k-b}^2.
\end{align*}
Since
\begin{align*}
L_{2k-b}^2=\sum_{\substack{p+q=b\\1\leq p\leq k\\1\leq q\leq
k}}\frac{(k+1-p)(k-p)q(q-1)}{(k+1)^2k^2}\sigma_p\sigma_q+\frac{b(b-1)}{(k+1)k}\sigma_b
\end{align*}
when $2\leq b\leq k+1$, and
\begin{align*}
L_{2k-b}^2=\sum_{\substack{p+q=b\\1\leq p\leq k\\1\leq q\leq
k}}\frac{(k+1-p)(k-p)q(q-1)}{(k+1)^2k^2}\sigma_p\sigma_q+\frac{(2k+2-b)(2k+1-b)}{(k+1)k}\sigma_{b-k-1}\sigma_{k+1}
\end{align*}
when $k+2\leq b\leq 2k$, then
\begin{equation}\label{e3.4}
L_{2k-b}=\left\{\begin{array}{ll}\dsum_{\substack{i+j=b\\1\leq i\leq
k\\1\leq j\leq
k}}\frac{(i-j+1)(k+1-i)j}{(k+1)^2k}\sigma_i\sigma_j&\\[12mm]
-\dfrac{b(b-1)}{(k+1)k}\sigma_b,&2\leq
b\leq k+1,\\[5mm]
\dsum_{\substack{i+j=b\\1\leq i\leq k\\1\leq j\leq
k}}\dfrac{(i-j+1)(k+1-i)j}{(k+1)^2k}\sigma_i\sigma_j&\\[12mm]
-\dfrac{(2k+2-b)(2k+1-b)}{(k+1)k}\sigma_{b-k-1}\sigma_{k+1},&
k+2\leq b\leq 2k.
\end{array}\right.
\end{equation}

\par
On the right hand side of \eqref{e3.3}, the coefficient of the term
$\alpha^{2k-b}(b\geq2)$ is
\begin{align*}
R_{2k-b}&=\sum_{\substack{i+j=b\\0\leq i\leq k\\0\leq j\leq
k}}\frac{(k+1-i)(k+1-j)}{(k+1)^2}\sigma_i\sigma_j-\sum_{\substack{p+q=b\\0\leq
p\leq k-1\\0\leq q\leq
k+1}}\frac{(k+1-p)(k-p)}{(k+1)k}\sigma_p\sigma_q\\
&:=R_{2k-b}^1-R_{2k-b}^2.
\end{align*}
\par
Obviously $R_{2k-b}=0$ when $b=0 (i=0,j=0,p=0,q=0)$. It is easy to
get $R_{2k-b}=0$ when $b=1 (i=1,j=0;i=0,j=1;p=1,q=0;p=0,q=1)$.
\par
When $2\leq b\leq k+1$,
\begin{align*}
R_{2k-b}^1=&\sum_{\substack{i+j=b\\1\leq i\leq k\\1\leq j\leq
k}}\frac{(k+1-i)(k+1-j)}{(k+1)^2}\sigma_i\sigma_j+\frac{2(k+1-b)}{k+1}\sigma_b,
\end{align*}
\begin{align*}
R_{2k-b}^2=&\sum_{\substack{p+q=b\\1\leq p\leq k\\1\leq q\leq
k}}\frac{(k+1-p)(k-p)}{(k+1)k}\sigma_p\sigma_q+\sigma_b+\frac{(k+1-b)(k-b)}{(k+1)k}\sigma_b,
\end{align*}
so
\begin{equation}\label{e3.5}
R_{2k-b}=\sum_{\substack{i+j=b\\1\leq i\leq k\\1\leq j\leq
k}}\frac{(k+1-i)(ik-jk+i)}{(k+1)^2k}\sigma_i\sigma_j-\frac{b(b-1)}{(k+1)k}\sigma_b,~~2\leq
b\leq k+1.
\end{equation}

\par
When $k+2\leq b\leq 2k$,
\begin{align*}
R_{2k-b}^1=&\sum_{\substack{i+j=b\\1\leq i\leq k\\1\leq j\leq
k}}\frac{(k+1-i)(k+1-j)}{(k+1)^2}\sigma_i\sigma_j,
\end{align*}
\begin{align*}
R_{2k-b}^2=&\sum_{\substack{p+q=b\\1\leq p\leq k\\1\leq q\leq
k}}\frac{(k+1-p)(k-p)}{(k+1)k}\sigma_p\sigma_q+\frac{(2k+2-b)(2k+1-b)}{(k+1)k}\sigma_{b-k-1}\sigma_{k+1},
\end{align*}
we get
\begin{align}\label{e3.6}
R_{2k-b}=&\sum_{\substack{i+j=b\\1\leq i\leq k\\1\leq j\leq
k}}\frac{(k+1-i)(ik-jk+i)}{(k+1)^2k}\sigma_i\sigma_j \\
&-\frac{(2k+2-b)(2k+1-b)}{(k+1)k}\sigma_{b-k-1}\sigma_{k+1},\quad
k+2\leq b\leq 2k. \nonumber
\end{align}

\par
Since
\begin{align*}
&(i-j+1)(k+1-i)j+(j-i+1)(k+1-j)i\\
=&(k+1-i)(ik-jk+i)+(k+1-j)(jk-ik+j),
\end{align*}
we have
\begin{equation}\label{e3.7}
\sum_{\substack{i+j=b\\1\leq i\leq k\\1\leq j\leq
k}}\frac{(i-j+1)(k+1-i)j}{(k+1)^2k}\sigma_i\sigma_j=\sum_{\substack{i+j=b\\1\leq
i\leq k\\1\leq j\leq
k}}\frac{(k+1-i)(ik-jk+i)}{(k+1)^2k}\sigma_i\sigma_j.
\end{equation}

\par
Using \eqref{e3.4},\eqref{e3.5},\eqref{e3.6} and \eqref{e3.7}, we
get equation \eqref{e3.3}.

\par
Now we prove \eqref{e3.1}. Using the equation(\ref{e3.2}) and
(\ref{e3.3}), we obtain
\begin{align*}
kn^2(S_{k;s}^2-S_{k-1;s}S_{k+1;s})(\bs
z)=&kn^2(E_k^2-E_{k-1}E_{k+1})(\bs z+\alpha\bs e)\\
=&\sum_{1\leq i<j\leq n}(z_i-z_j)^2\prod_{l\neq
i,j}(z_l+\alpha)^2\geq0.
\end{align*}
Specifically, when $s=3, k=4, n=5$, we have
\begin{align*}
&100[S_{4;3}^{2}(\bs z)- S_{3;3}(\bs z)S_{5;3}(\bs z)]\\
=&(z_1-z_2)^2(z_3+\alpha)^2(z_4+\alpha)^2(z_5+\alpha)^2+(z_1-z_3)^2(z_2+\alpha)^2(z_4+\alpha)^2(z_5+\alpha)^2\\
&+(z_1-z_4)^2(z_2+\alpha)^2(z_3+\alpha)^2(z_5+\alpha)^2+(z_1-z_5)^2(z_2+\alpha)^2(z_3+\alpha)^2(z_4+\alpha)^2\\
&+(z_2-z_3)^2(z_1+\alpha)^2(z_4+\alpha)^2(z_5+\alpha)^2+(z_2-z_4)^2(z_1+\alpha)^2(z_3+\alpha)^2(z_5+\alpha)^2\\
&+(z_2-z_5)^2(z_1+\alpha)^2(z_3+\alpha)^2(z_4+\alpha)^2+(z_3-z_4)^2(z_1+\alpha)^2(z_2+\alpha)^2(z_5+\alpha)^2\\
&+(z_3-z_5)^2(z_1+\alpha)^2(z_2+\alpha)^2(z_4+\alpha)^2+(z_4-z_5)^2(z_1+\alpha)^2(z_2+\alpha)^2(z_3+\alpha)^2\geq0.
\end{align*}

It is obvious to see that the equality holds when
$z_{1}=z_{2}=\cdots=z_{n}$, or any two elements of
$z_{1},z_{2},\cdots,z_{n}$ are equal to $-\alpha$.

\par
We assume that
 $z_{1}=z_{2}=-\alpha$, note that $n=k+1$, then by \eqref{e2.6},
\begin{align*}
&C_{k+1}^1E_1(\bs z)=\sigma_1(\bs z)=-2\alpha+\sigma_1({\bs z}'),
\\
&C_{k+1}^iE_i(\bs z)=\sigma_i(\bs z)=\alpha^2\sigma_{i-2}({\bs
z}')-2\alpha\sigma_{i-1}({\bs z}')+\sigma_i({\bs z}'),~~~
\hbox{for}~~~ 2\leq i\leq k-1,\\
&C_{k+1}^kE_k(\bs z)=\sigma_k(\bs z)=\al^2\sigma_{k-2}({\bs
z}')-2\alpha\sigma_{k-1}({\bs z}'),\\
&C_{k+1}^{k+1}E_{k+1}(\bs z)=\sigma_{k+1}(\bs
z)=\al^2\sigma_{k-1}({\bs z}').
\end{align*}
Note again that $n=k+1, s=k-1$, using the identities above and the
define of $S_{k;s}$ in \eqref{e1.7}, we have
\begin{align*}
S_{k;s}(\bs z)=&C^{0}_{s}E_{k}(\bs x)+C^{1}_{s}\alpha E_{k-1}(\bs
x)+\cdots+C^{s}_{s}\alpha^s E_{k-s}(\bs x)\\
=&\dfrac{C_{k-1}^0}{C_{k+1}^k}[\alpha^2\sigma_{k-2}({\bs
z}')-2\alpha\sigma_{k-1}({\bs
z}')]\\
&+\dfrac{C_{k-1}^1}{C_{k+1}^{k-1}}\alpha[\alpha^2\sigma_{k-3}({\bs
z}')-2\alpha\sigma_{k-2}({\bs z}')+\sigma_{k-1}({\bs z}')]\\
&+\cdots+\dfrac{C_{k-1}^{k-1}}{C_{k+1}^1}\alpha^{k-1}[-2\alpha+\sigma_1({\bs
z}')]\\
=&\dfrac{-2\alpha}{k(k+1)}[\sigma_{k-1}({\bs
z}')+\alpha\sigma_{k-2}({\bs z}')+\cdots+\alpha^{k-1}].
\end{align*}
Similarly, we have
\begin{align*}
S_{k-1;s}(\bs z)=&\dfrac{2}{k(k+1)}[\sigma_{k-1}({\bs
z}')+\alpha\sigma_{k-2}({\bs z}')+\cdots+\alpha^{k-1}],\\
S_{k+1;s}(\bs z)=&\dfrac{2\alpha^2}{k(k+1)}[\sigma_{k-1}({\bs
z}')+\alpha\sigma_{k-2}({\bs z}')+\cdots+\alpha^{k-1}].
\end{align*}
So
\begin{equation*}
  \frac{S_{k;s}(\bs z)}{S_{k-1;s}(\bs z)}=\frac{S_{k+1;s}(\bs z )}{S_{k;s}(\bs z)}=-\al.
\end{equation*}
We completed the proof of Lemma \ref{lem3.1}.

\par
{\bf Proof of Theorem \ref{th5}}  Combine Lemma~\ref{lem2.2} with
Lemma~\ref{lem3.1}, and use the method in the proof of \eqref{e1.3}
of the Theorem~\ref{th1}, we can similarly prove \eqref{e1.8} of the
Theorem~\ref{th5}.

\par
The proof of \eqref{e1.10} is similar to \eqref{e1.5}, where the key
step is to prove
\begin{equation}\label{e3.8}
S_{k;s}^{2}(\bs x)- S_{k-1;s}(\bs x)S_{k+1;s}(\bs x)\geq 0, \quad
k=1,\cdots, s,
\end{equation}
with the condition $\alpha \geq 0, E_1(\bs x)\geq 0, E_2(\bs x)\geq
0,\cdots,E_s(\bs x)\geq 0$.
\par
Now we prove \eqref{e3.8}. Let $1\leq k<s$, by the define of
$S_{k;s}$ in \eqref{e1.7},
\begin{align*}
S_{k;s}(\bs z)&=\sum_{i=0}^kC^{i}_{s}\alpha^{i}E_{k-i}(\bs z),\\
S_{k-1;s}(\bs z)&=\sum_{i=0}^{k-1}C^{i}_{s}\alpha^{i}E_{k-1-i}(\bs z),\\
S_{k+1;s}(\bs z)&=\sum_{i=0}^{k+1}C^{i}_{s}\alpha^{i}E_{k+1-i}(\bs
z).
\end{align*}
From the expansion above, we only need to show that all the terms
containing $\alpha$ at left hand sides of \eqref{e3.8} are
nonnegative. The coefficient of the term $\alpha^{b}$ is
\begin{align*}
L_{b}&=\sum_{\substack{i+j=b\\0\leq i\leq k\\0\leq j\leq
k}}C_s^{i}C_s^{j}E_{k-i}E_{k-j}-\sum_{\substack{p+q=b\\0\leq p\leq
k-1\\0\leq q\leq k+1}}C_s^{p}C_s^{q}E_{k-1-p}E_{k+1-q}\\
&:=L_{b}^1-L_{b}^2.
\end{align*}
\par
Obviously, from the expression of $L_{b}$ above, we have
$$L_0=E_k^2-E_{k-1}E_{k+1}\geq 0,\quad
L_{2k}=C_s^kC_s^k-C_s^{k-1}C_s^{k+1}\geq 0.
$$
and
\begin{align*}
L_k=&\sum_{i=0}^kC_s^{i}C_s^{k-i}E_{k-i}E_{i}-\sum_{i=0}^{k-1}C_s^{i}C_s^{k-i}E_{k-i-1}E_{i+1}\\
\geq&\sum_{i=0}^{k-1}(C_s^{i}C_s^{k-i}E_{k-i}E_{i}-C_s^{i}C_s^{k-i}E_{k-i-1}E_{i+1})\\
=&\sum_{i=0}^{[k]/2}(C_s^{i}C_s^{k-i}E_{k-i}E_{i}+C_s^{k-1-i}C_s^{i+1}E_{i+1}E_{k-1-i}\\
&-C_s^{i}C_s^{k-i}E_{k-i-1}E_{i+1}-C_s^{k-1-i}C_s^{i+1}E_{i}E_{k-i})\\
:=&\sum_{i=0}^{[k]/2}R_i.
\end{align*}
In the second equal sign above, we used a summation method from both
ends to the middle for $i$. Since $i\leq [k]/2$, which implies that
$k-i\geq i$. If $k-i> i$,
$$
E_{i+1}E_{k-1-i}\geq E_{i}E_{k-i},\quad C_s^{k-1-i}C_s^{i+1}\geq
C_s^{i}C_s^{k-i},
$$
we get $R_i\geq 0$.  If $k-i=i$,
$$
R_i=C_s^iC_s^iE_i^2+C_s^{i-1}C_s^{i+1}E_{i+1}E_{i-1}-C_s^iC_s^iE_{i-1}E_{i+1}-C_s^{i-1}C_s^{i+1}E_i^2\geq
0.
$$
So we proved $L_k\geq 0$.

\par
When $0<b<k$,
\begin{align*}
L_{b}^1&=\sum_{i=0}^bC_s^{i}C_s^{b-i}E_{k-i}E_{k-b+i},\\
L_{b}^2&=\sum_{i=0}^bC_s^{i}C_s^{b-i}E_{k-i-1}E_{k-b+i+1}.
\end{align*}
Using the similar method to prove that $L_k\geq 0$, we can get
$$
L_b=L_{b}^1-L_{b}^2\geq 0.
$$

\par
When $k<b<2k$,
\begin{align*}
L_{b}^1&=\sum_{i=b-k}^k C_s^{i}C_s^{b-i}E_{k-i}E_{k-b+i},\\
L_{b}^2&=\sum_{p=b-(k+1)}^{k-1}
C_s^{p}C_s^{b-p}E_{k-p-1}E_{k-b+p+1}=\sum_{i=b-k}^{k}
C_s^{i-1}C_s^{b-i+1}E_{k-i}E_{k-b+i}.
\end{align*}
The proof of $L_b\geq 0$ is also similar as $L_k$.
\par
Thus we obtain \eqref{e3.8},  and the Theorem \ref{th5} is proved.

\par

\end{document}